\documentclass[titlepage,11pt,leqno]{article}
\usepackage{latexsym}
\usepackage{amsmath}
\usepackage{amsfonts}
\usepackage{threeparttable}
\usepackage[dvips]{graphicx}
\def\lel{\label}
\def\r{\ref}

\begin{document}
\baselineskip 16.5pt
\bibliographystyle{unsrt}
\newcommand {\apgt} {\ {\raise-.5ex\hbox{$\buildrel>\over\sim$}}\ }
\newcommand {\aplt} {\ {\raise-.5ex\hbox{$\buildrel<\over\sim$}}\ }
\newtheorem{theorem}{\indent \sc Theorem}
\newtheorem{corollary}{\indent \sc Corollary}
\newtheorem{lemma}{\indent \sc Lemma}
\newtheorem{proposition}{\indent \sc Proposition}
\newtheorem{remark}{\indent \sc Remark}
\newcommand{\phif}{\textsc{igf}}
\newcommand{\phgf}{\textsc{gf}}
\newcommand{\fix}{$\textsc{gf}_0$}
\newcommand{\mb}[1]{\mbox{\bf #1}}
\newcommand{\Exp}[1]{\mbox{E}\left[#1\right]}
\newcommand{\Var}[1]{\mbox{{\rm var}}\left[#1\right]}
\newcommand{\Cov}[2]{\mbox{{\rm cov}}\left[#1,\;#2\right]}
\newcommand{\pr}[1]{\mbox{P}\left[#1\right]}
\newcommand{\pp}[0]{\mathbb{P}}
\newcommand{\ggg}[0]{\mathbb{G}}
\newcommand{\ee}[0]{\mbox{E}}
\newcommand{\re}[0]{\mathbb{R}}
\newcommand{\argmax}[0]{\mbox{argmax}}
\newcommand{\argmin}[0]{\mbox{argmin}}
\newcommand{\ind}[0]{\mbox{\Large\bf 1}}
\newcommand{\narrow}{\stackrel{n\rightarrow\infty}{\longrightarrow}}
\newcommand{\ol}[1]{\overline{#1}}
\newcommand{\avgse}[1] { \bar{\hat{\sigma}}_{#1} }
\newcommand{\mcse}[1]  { \sigma^{*}_{#1} }
\newcommand{\po}{\textsc{po}}
\newcommand{\ph}{\textsc{ph}}
\noindent
\noindent \begin{center}{\LARGE\bf Model Selection and Multiple Testing - A \textbf{Bayes} and Empirical Bayes Overview and some New Results}
\vspace{0.1in}

\noindent   Ritabrata Dutta$\mbox{}^{a}$, Malgorzata Bogdan$\mbox{}^{a}$$\mbox{}^{b}$ and Jayanta K Ghosh$\mbox{}^{a}$$\mbox{}^{c}$

\noindent $\mbox{}^{a}$ {\it Department of Statistics, Purdue University}

\noindent $\mbox{}^{b}$ {\it Institute of Mathematics and Computer Science, Wroclaw University of
Technology, Poland}

\noindent $\mbox{}^{c}$ {\it Indian Statistical Institute, India.}

\end{center}
\vskip .2in
{\small
\noindent{\bf Abstract}:
We provide a brief overview of both Bayes and classical model selection.
We argue tentatively that model selection has at least two major goals,
that of finding the correct model or predicting well, and that in general
both these goals may not be achieved in an optimum manner by a single model
selection rule. We discuss, briefly but critically, through a study of
well-known model selection rules like AIC, BIC, DIC and Lasso, how these
different goals are pursued in each paradigm. We introduce some new
definitions of consistency, results and conjectures about consistency in high
dimensional model selection problems. Finally we discuss some new or recent
results in Full Bayes and Empirical Bayes multiple testing, and cross-validation.
We show that when the number of parameters tends to infinity at a smaller rate
than sample size, then it is best from the point of view of consistency to use
most of the data for inference and only a negligible proportion to make an
improper prior proper.
}

\vspace{0.2in}

\vspace{0.05in}

\section{Introduction}

We provide a brief history of model selection starting with the work of Jeffreys (\textbf{1961}), 
Cox (1961, 1962), Akaike (1974) and \textbf{Schwarz} (1978) at the begining of Section \ref{s2}. 
Each of these writers had a well-defined purpose for model selection, which was either 
choosing the true model or choose a model that predicts optimally in some well-defined sense. 
The actual method adopted in each case also depends on the paradigm chosen, Bayesian or classical.

For a given $\theta$, if the goal is to choose the true model, one hopes to be consistent 
in the sense that one chooses the model containing true $\theta$. The precise definition 
is given later in terms of the Bayes factor or likelihood ratio. Similarly, for optimal 
prediction, one needs to develop an oracle, i.e., a lower bound to risk of all candidate 
procedures, which is asymptotically attained by the model selection rule under consideration. 
These ideas are illustrated by AIC and BIC. Then we compare AIC and BIC, and also discuss more 
recent rules like DIC and Lasso. These discussions are in Subsections \ref{s21}, \ref{s22} and \ref{s23}.

Subsection \ref{s24} deals briefly with Parametric Empirical Bayes model selection as
introduced by George and Foster (2000). Their formulation takes care of issues
relating to both complexity and multiplicity. We also discuss some optimality
results of Mukhopadhyay and Ghosh (2003) in the setting of George and Foster (2000).

In Section \r{s3}, we discuss high dimensional Bayes model selection. We review what
is known, propose new definitions of consistency, which may be easier to prove
and which make better sense in situations such as Scott and Berger (2010). We
also prove new results including a theorem under the new definition of
consistency and make a couple of conjectures on consistency as usually defined.

Section \r{s4} deals with Bayes prediction and cross-validation,
with more stress on the latter. Our main result in this section is Theorem 2,
which is the the result mentioned in the Abstract. We also argue that the recently
proposed cross-validatory Bayes factor of Draper and \textbf{Krnjajic} (2010) is actually
not very different from the usual Bayes factor for fixed $p$ and moderate or
large $n$, but leads to interesting differences in model selection. We point this
out in the hope that these major issues will be thoroughly discussed among Bayesians.
Section 5 introduces the reader to the closely related topic of multiple testing,
which can in fact be regarded as model selection, and has been one of the fastest
growing areas of theoretical and applied statistics.

Each section has something new to offer, a history and an overview in Sections \textbf{2} and 5,
an attempt to find new directions in Section 3 and a new result in Section 4.

\section{Early History of Bayes and Classical Model Selection}\label{s2}

To the extent that Bayes model selection, with the goal of choosing the correct
 model, is an extension of or identical with testing of two hypotheses, Bayes model
selection may be said to have begun with the work of Jeffreys in 1939, see
Jeffreys (1961). The Bayes test, as well as an approximation to it due to Lindley,
is mentioned in Cox (1961). Cox (1961, 1962) pioneered model selection in classical
statistics, but without bringing in the novel concept of penalization for complexity of
a model.

Both in Cox (1961) and Cox (1962), one considers two separated hypotheses, i.e.
hypotheses or models having the property that no density in one is obtainable as
a limit of densities in the other. Cox does not specify the notion of convergence
but in a follow up study, Ghosh and Subramanyam (1975), suggest that the limit
may be taken in the sense of convergence in $L_1$-norm. This is equivalent to
requiring the two sets of densities can be covered by two disjoint $L_1$-open
sets. An alternative definition, depending on n, is given towards the end of the previous reference.

Essentially, Cox's model selection rule is based on the maximized likelihood of
data under each model. Subramanyam and Ghosh, vide Subramanyam (1979), show that
the true model is rejected with exponentially small probability. To show this, one
has to use results on large deviations. They verify the conditions for one of Cox's
examples, where the true model is either Geometric or Poisson.

It is remarkable that both in Jeffreys (1961) and Cox (1961, 1962) the goal is the
same - to choose the correct model. It is equally remarkable that, starting with AIC,
choosing the correct model has not been a goal in classical model selection, at least
not explicitly.

If one reads Akaike's early papers and resolves the ambiguities in the light of the
subsequent pioneering theoretical papers of Shibata (1984), it becomes clear that
the objective is to choose a model so that one predicts optimally the data on the
dependent variables in an exact new replicate of the design for the given data.
Assuming the true model is more complex than the assumed linear models, Shibata's
calculations, and also the calculations in Li (1987) and Shao (1997), show that
the special form of penalty associated with AIC has an intimate role in this kind
of optimal prediction. Optimality is proved through a lower bound to the predictive
loss or risk of all so called penalized likelihood rules that may be used in this
problem and then showing AIC attains the lower bound asymptotically because of its
special penalty. This lower bound would now be called an oracle, a term which was
probably introduced later formally in the model selection literature by
Johnstone and Donoho (1994). It has been shown by van de Geer (2005) that
there is a close connection between the penalty of a predictive model selection
criteria and the oracle it may attain.

In Bayes model selection, the picture is far more mixed. One sees papers of both
kinds, selecting the correct model or predicting well. The rich literature on Bayes
model selection is reviewed in Clyde and George (2004).

A historian of model selection may well speculate about this curious divergence
in development. To us it seems this may be due to the clarity of all aspects,
including goals of inference, that one finds readily in the Bayesian paradigm.

In the first of the next two subsections, we first introduce Bayes model
selection through Bayes factors and then \textbf{Schwarz}'s BIC. In the next subsection
we discuss a few popular model selection rules based on penalized likelihood.

\subsection{Bayes Rule for Selection of True Model: BIC, Consistency
in Bayes model selection}\label{s21}
Assume first that we have two models $M_i,\ i=1,2 $. With $M_i$, we
associate a parameter space $\Theta_i$, a density
$f_i(x_1,x_2,\cdots,x_n|\theta_i)$ for the data $\mathbf{x}=(x_1,x_2,\cdots,x_n\mathbf{)}\in
\mathbb{R}_n$ and a prior $\pi_i(\theta_i)$ for $\theta_i$. We assume the
families $\{f_i(\mathbf{x}|\theta_i),\theta_i\in \Theta_i\},i =1,2$ are separated
in the sense explained earlier.

Assume the Bayesian has also prior probabilities $\lambda_1,\lambda_2$ for $M_1$
and $M_2$, $\lambda_1+\lambda_2=1$. Let
\begin{eqnarray*}
m_i(\mathbf{x})&=&\mbox{marginal density of}\ \mathbf{x} \ \mbox{under } M_i
=\int_{\Theta_i}f_i(\mathbf{x}|\theta_i)\pi_i(\theta_i)d\theta_i.
\end{eqnarray*}
Then the posterior probability of $M_i$, given $\mathbf{x}$, is proportional
to $\lambda_i m_i(\mathbf{x})$. Hence one would select $M_2$ if it is more likely
than $M_1$ given data, i.e., if
\begin{eqnarray*}\label{1}
\frac{\lambda_2 m_2(\mathbf{x})}{\lambda_1 m_1(\mathbf{x})}>1,
\end{eqnarray*}
and the other way if the ratio is $<1$. Jeffreys (1961) has suggested a scale of
numerical values for different levels of relative credibility of the two models.
The Bayes rule is easily extended when we have k separated models.

The usual default choice for $\lambda_1,\lambda_2$ is
$\lambda_1=\lambda_2=\frac{1}{2}$, in which case our selection criterion becomes:
\begin{eqnarray*}
\mbox{Select}\ M_2, \ \mbox{if} \ \frac{m_2}{m_1}>1.
\end{eqnarray*}
The ratio $\frac{m_2}{m_1}$ is called the Bayes factor and denoted by $BF_{21}$.

We usually assume that each model specifies iid models
$f_i(\mathbf{x}|\theta_i)=\prod_1^nf_i(x_j|\theta_i)$.

While Bayes model selection is straightforward in principle, numerical calculation
of the Bayes factor is not. One needs Reversible Jump MCMC or Path Sampling,
Andrieu et. al. (2004), Dutta and Ghosh (2011). \textbf{Schwarz}'s (1978) BIC
provides a convenient approximation, assuming fixed dimension $p_i$ of
$\Theta_i$ and suitable regularity conditions, and letting sample size $n \to \infty$,
\begin{eqnarray*}
\log{m_i}=\sum^n_{\mathbf{j=1}}\log{f_{i}(x_j|\hat{\theta}_i)}-\frac{p_i}{2}\log{n}+O_p(1),
\end{eqnarray*}
where $p_i$, equal to the dimension of $\Theta_i$, is a measure of the
complexity of $M_i$, and $\hat \theta_i$ is the MLE under model $M_i$. Thus $\log{m_i}$ can be approximated by
\begin{eqnarray}\label{4}
BIC(M_i)=\sum^n_{\mathbf{j=1}}\log{f_i(x_j|\hat{\theta}_i)}-\frac{p_i}{2}\log{n},
\end{eqnarray}
which may be interpreted as a penalized maximum likelihood corresponding to a model.
The bigger the dimension of $\Theta_i$ the bigger is the penalty, in tune with the
scientific principle of parsimony. One selects a model by maximizing equation (\ref{4})
with respect the model $M_i$, i.e., by choosing suitable  density  $f_i$.
This is the user's guess for the unknown true model.

Thus Bayesian model selection through Bayes factors automatically obeys the
principle of parsimony. It tries to compensate for the fact that the bigger the $\Theta_i$, the bigger we expect the maximum likelihood to be.

As theoretical validation of use of Bayes factor in a problem, one tries to prove consistency in some sense. Below is the usual definition adopted by Bayesians, but we do not know of a reference. We first consider the case of two separated models $M_1$, $M_2$ with our usual notation. The nested case will be discussed in the next section.

\textbf{Definition of Consistency}. We say consistency holds at
$\theta\in \Theta_2$ if under $\theta$, $\log{BF_{21}}\to \infty$,
in probability. On the other hand consistency holds at $\theta\in \Theta_1$,
if $\log{BF_{21}}\to -\infty$, in probability.

Consistency in this sense will appear in our discussions many times.
Also in Section 3 a definition of an alternative notion will be given. It
will be used only in the example following that definition.

\subsection{AIC and BIC}\label{s22}

We start this subsection by digressing a little on nested models and some
issues that arise when one of two nested models have to be selected. $M_1$ is
said to be nested in $M_2$ if $\Theta_1\subset\Theta_2$ and for $ \theta\in \Theta_1$,
$f_1(x_j|\theta)=f_2(x_j|\theta)$. This is a very different situation from the
separated models that we have been considering earlier. If $\theta\in \Theta_2$ but
$\notin \Theta_1$, we may say $M_2$ is true and the density is $\mathbf{f_2(x|\theta)}$. But
what about the case where $\theta\in \Theta_1$ and hence $\theta\in \Theta_2$ also?
Here both models are correct, so which one do we choose? Bayesians suggest that on
grounds of parsimony one should choose the smaller model $M_1$.

What about priors over $M_1$ and $M_2$? Typically, we have priors $\pi_1(\theta_1)$
on $M_1$ and $\pi_2(\theta_2)$ on $M_2$ such that $\pi_2(\theta_2)$ is a density
with respect to Lebesgue measure on $\Theta_2$ and hence assigns zero probability
to $\Theta_1$. Typically, $\pi_1(\theta_1)$ has a density with respect to Lebesgue
measure on $\Theta_1\subset\Theta_2$, so it assigns zero probability to
$\Theta_2\cap\Theta_1^c$.

The problem of selection of one of two nested models has puzzled philosophers of
science. They feel one should always choose the bigger model $M_2$ as true.
How can $M_2$ be rejected if $M_1$ is true? For their views, as well as what
Bayesians have to say in response about logical consistency of nested model
selection, see Chakrabarti and Ghosh (2011). They try to explain through Galileo's
famous experiment at Pisa. The smaller model represents Gallileo's view, while the
bigger model represents the prevalent view at that time about falling bodies.
The two priors represent idealized versions of two such views. Typically when we
deal with linear models, choosing means or choosing non-zero regression coefficients
in regression problems, we have to select from among models, some of which are
nested in others. This is the most common model selection problem.

We now turn to a comparison of AIC and BIC. Historically BIC (Bayes Information
Criterion) is the second penalized maximum likelihood rule. The first is
AIC (Akaike Information Criteria). For linear models with normal error and
known variance $\sigma^2$, AIC proposed by Akaike (1974) is
\begin{eqnarray}\label{aic}
AIC(M_i)=2\sum^n_{\mathbf{j=1}}\log{f_i(x_j|\hat{\theta}_i)}-2p_i.
\end{eqnarray}
One has to \textbf{multiply} BIC by two or divide AIC by two to make them comparable.
The expression (\ref{aic}) is maximized over models to get the best predictive model.
Prediction is made by using the mle under the chosen model.

The first term of AIC is taken to be twice the more intuitive first term of BIC,
because the difference of the first terms for AIC for two models has
asymptotically a $\chi^2$-distribution if the smaller model is true.
We could have done the same for BIC by doubling the penalty.

We recall that though AIC and BIC look similar, they are meant to do very
different things. AIC predicts optimally, while, as pointed out above, BIC
is a Bayesian criterion for selecting the true model. Each will perform poorly if
used for a purpose for which it was not meant.

We now turn to an example of a simple linear model, with normal error N(0,1)
\begin{eqnarray*}
Y_i=\mu+\epsilon_i,\quad i=1,2,\cdots, n.
\end{eqnarray*}
Suppose we wish to choose between $M_1:\mu=0$ and $M_2:\mu$ is arbitrary, i.e.,
$\Theta_1=\{0\}$, $\Theta_2=\mathbb{R}$. If we use AIC to select the true model,
then it is easy to verify that we are actually using a test with Type 1 error
= $P\{\chi^2_1>2\}>.05$, i.e., the  test is more liberal than the usual most liberal
test. So, there is no parsimony if we use AIC in a testing problem.

On the other hand suppose we use BIC. BIC is very parsimonious and can be shown to
choose the correct model with probability tending to one, at least for fixed $p$
and $n \to \infty$. If we have linear models and our goal is prediction, AIC will
be predictively optimal as discussed earlier. For ease of reference we consider
Shao (1997). Note that Shao's theorem is proved in the context of general linear
models, and the main assumption is that the true linear model is not in the model
space. Shao's theorem is an illustration of Box's famous remark that all models are
false but some are useful. Moreover Shao's theorem shows AIC helps one identify the
most useful model from the point of view of prediction. The proof is based on an
oracle. Chakrabarti and Ghosh (2007) have used a similar oracle property of AIC
to show heuristically that it is an adaptive, asymptotically minimax, estimate of
the unknown non-parametric regression function $f$ in the following model,
\begin{eqnarray*}
Y_i=f(t_i)+\epsilon_i,\ \epsilon_i \stackrel{\mathrm{iid}}{\sim}  N(0,\sigma^2),
\,\,\ i=1,\cdots, n, \quad t_i=\frac{i}{n}.
\end{eqnarray*}
Computation in Chakrabarti (2004) shows AIC is competetive with popular best methods in nonparametric regression.

On the basis of comparison of AIC and BIC, we suggest tentatively that model
selection rules should be used for the purpose for which they were introduced.
If they are used for other problems, a fresh justification is desirable. In one
case, justification may take the form of a consistency theorem, in the other some
sort of oracle inequality. Both may be hard to prove. Then one should have
substantial numerical assessment over many different examples.

\subsection{DIC and Lasso}\label{s23}

In this subsection we make a few remarks about DIC and Lasso, both of which have
been relative late comers, but both of which are very popular.

Deviance Information Criteria (DIC) has been developed by Bayesian authors, (Spiegelhalter et. al., 2002), in an effort to generalize AIC, using the Kullback-Liebler divergence instead of squared error loss. To some extent this has been done earlier also. However, DIC also tries to find a form of penalty that can take care of \textbf{hierarchical} models as well as latent parameter models. It seems DIC is really suitable for prediction but it has also been used for choosing the correct model. Also, there have been doubts about the penalty. Interested readers should read both the discussion and the reply of the authors.

For variable selection in linear regression problems shrinkage estimators
like ridge regression are very popular. Following the idea of shrinkage as
in ridge regression, the Lasso method introduced by Tibshirani(1994), tries
to minimize the least square error of the regression with an upper bound on
the $L_1$ norm of the parameter vector.  The estimate is defined by,
\begin{eqnarray*}
\widehat{\beta}^{lasso}&=&\argmin_{\beta}\sum^N_{i=1}(y_i-\beta_1-\sum^p_{j=2}x_{ij}
\beta_j)^2\, , \qquad
\mbox{subject to}\, \, \sum^p_{j=2}|\beta_j|\leq t\,.
\end{eqnarray*}
In this connection it is worth mentioning (Candes and Tao, 2007) that according
to unpublished work by Donoho, for most noiseless large undetermined systems of
linear equations the minimal $l_1$-norm solution is also the sparsest solution.
Donoho's work provides deep insight about the success of Lasso. The upper bound $t$
for the $L_1$-penalty of the parameter-vector controls the amount of shrinkage.
Lasso chooses subsets of variables depending on the tuning-parameter $t$. When $t$
is very small, then almost all the parameters are zero; similarly for large
enough $t$, the parameter estimate $\widehat{\beta}$ is the same as the least square
estimate. The shrinkage constraint makes the solution nonlinear in $y_i$, needing a
quadratic programming algorithm to compute the estimate.  Efron et. al.(2004)
have shown that Lasso is closely related to another novel shrinkage estimation
scheme LARS, introduced by them. In a general setup, they show that the subset
selected by LARS and Lasso are similar. Here we will concentrate on the most
popular of them, namely, Lasso, its solution path obtained from the LARS
algorithm for fast variable selection. This has become the standard method for
Lasso. The tuning parameter $t$ is chosen by minimizing the cross-validation error.
To the extent that new methods like Lasso have become popular, stepwise all
subsets model selection has become much less used.

Recently Bickel et. al. (2009) have derived an Oracle property for Lasso under a
sparsity assumption and some other rather stringent conditions. But we get a lot of
insight about Lasso from its oracle inequality, as pointed out towards the end of this
paragraph. The Lasso constraint $\sum{|\beta_j|}\leq t$ is equivalent to the addition
 of a penalty term $r\sum|\beta_j|$ to the residual sum of squares (Murray et. al. (1981)). 
While an explicit 
mathematical relation between $t$ and $r$ is not available, the basic idea of convex 
optimization makes it easy to move from the one to the other. One may use both versions of 
Lasso. This can be written as follows,
\begin{eqnarray*}
\widehat{\beta}^{lasso}=\argmin_{\beta}\sum^N_{i=1}(y_i-\beta_1-\sum^p_{j=2}x_{ij}
\beta_j)^2+r\sum^p_{j=2}|\beta_j|\;\;.
\end{eqnarray*}
According to Bickel et. al.(2009), when the errors $\epsilon_i$ are
independent $N(0,\sigma^2)$ random variables with $\sigma^2>0$, and all the
diagonal elements of the {\bf Gram} matrix $X'X/n$ are equal to 1, then under
some additional conditions on the Gram matrix,  $r=A\sigma\sqrt{\frac{\log(p)}{n}}$
and $A>2\sqrt{2}$, with probability $1-p^{1-\frac{A^2}{8}}$, we have
\begin{eqnarray*}
|\widehat{\beta}^{lasso}-\beta_0|_1&\leq &\frac{16A}{c(s)}\sigma 
s \sqrt{\frac{\log(p)}{n}}\;,\\
|X(\widehat{\beta}^{lasso}-\beta_0)|^2_2&\leq &\frac{16A^2}{c(s)}\sigma^2 s 
\log(p)\;,\qquad 
\#\{\widehat{\beta}_j^{lasso}\neq 0\}\leq \frac{64}{c_1(s)}s\;,
\end{eqnarray*}
where $s$ is the number of non-zero components in $\beta_0$, and $c(s)$, $c_1(s)$ are constants 
depending on $s$ and the Gram matrix. They have very similar oracles for the Dantzig selector of 
Candes and Tao (2007), suggesting both methods achieve similar goals. The oracle also makes clear 
that the penalty should change with sparsity.

The penalized model selection schemes like Lasso have also been implemented in the Bayesian set 
up by Park and Casella (2008) and Kyung et. al. (2010). They have shown that with the proper 
choice of prior distribution, the posterior distribution introduces a penalty term, e.g. Laplace 
(double-exponential) distribution as prior associates a penalty term that is same as Lasso, thus 
replicating penalized model selection schemes by Bayesian methods. They also demonstrated that the 
generalized Lasso schemes including Fused Lasso, Elastic net, Group Lasso and others can be 
implemented in Bayesian set up by proper choice of a prior in Bayesian Lasso. The discussions 
about the consistency, standard error, performance and comparison of Bayesian Lasso with 
standard penalized schemes are also illuminating (Kyung et. al., 2010).

While on the subject of consistency in the context of linear models, we like to mention a few 
papers chosen from the emerging literature on consistency used in the sense of approximation 
from a given dictionary of functions. Much of the new vocabulary has come from Machine Learning, 
but almost all the papers we cite have appeared in statistical journals, some have a Bayesian 
flavor, Bunea et. al. (2007), Bunea(2008), Zhao and Yu (2007), Zou(2006). Bunea (2008) contains 
many references on sparsity, oracles, information theoretic limits.Several of these relate to the 
Lasso.

\subsection{Predictive Bayesian Model Selection - Model average and median model}
\label{s24}

If the object of Bayes model selection is not to select the true model but 
predict future observations well, the utility or loss changes dramatically
 from 0-1 loss to a conditional expectation of squared error loss of predictor 
of future observations $x_f$. Instead of the squared error one may choose other 
suitable loss functions for prediction. The Bayesian solution of this problem is 
straightforward but very different from the posterior mode selecting models, which 
is what we have discussed so far.

For prediction it is best not to select a single model but rather average
over best predictions from all models, i.e., use the so called predictor based on
model average. Assuming conditional independence of current data $x$ and future data
$x_f$, given $\theta$, calculate the model average prediction
\begin{eqnarray*}
\frac{\sum{P(M_j|X_j)\int{E(x_f|\theta_j,M_j)p(\theta_j|M_j)d\theta_j}}}
{\sum{P(M_j|X_j)}}\, .
\end{eqnarray*}

If, on the other hand, one wants to choose a suitable model $M_j$ and then
use the predictor assuming $M_j$ is true, then it turns out that the median
model of Barbieri and Berger (2004) is optimal under orthogonality assumptions.
A quick introduction to both model average and median model is available in
Ghosh et. al. (2006, ch 9). As mentioned earlier, a posterior quantile model
selection rule and its optimality is studied in Mukhopadhyay and Ghosh (2003).

Finally computation is a major concern, specially when both the dimensions of
models and the number of models are very large. Many algorthims are available,
see e.g. Clyde and George (2004).

\subsection{Parametric Empirical Bayes (PEB) Model Selection}\lel{s25}

A great advantage of the PEB approach is that it often makes theoretical
calculations surprisingly easy. Such calculations would be almost impossible
with a Full Bayes or Classical approach. Typically, PEB is applied when the
likelihood is a product of the likelihood of individual parameters. Typically,
this structure is utilized by choosing a prior that makes the parameters independent,
given some hyperparameters, through which borrowing of strength across different
$X_i$'s occur. Up to this, PEB and Hierarchical Bayes (HB) have the same structure.
Now they part ways. HB puts a prior on the hyperparameters while PEB treats them
as unknown hyperparameters to be estimated as follows. First integrate out the lower level parameters. The resulting likelihood for hyperparameters is then
maximized to get what are called Type 2 MLE's. Alternatively one gets Best
Unbiased Estimators as in Morris (1983). The relative simplicity of PEB comes
at a price. There is \textbf{no} adjustment for the uncertainty in the inference due to
 hyperparameters. HB takes care of the full uncertainty.

We \textbf{illustrate} with the PEB approach to model selection due to George and Foster (2000),
as applied to nested models, and some optimality results for nested models due to
Mukhopadhyay and Ghosh (2003) that follow from the approach of George and Foster
(2000). Throughout, the covariates are assumed to be mutually orthogonal.
The data consist of independent r.v.'s $Y_{ij},\ i=1,2,\cdots,p,\ j=1,2,\cdots,r$.
There are $p$ models $M_q,\ 1\leq q \leq p$. Under $M_q$,
\begin{eqnarray*}
Y_{ij}&=&\beta_i+\epsilon_{ij},\hskip .2in 1\leq i \leq q,\ j=1,2,\cdots,r\\
&=&\epsilon_{ij},\hskip .3in q+1 \leq i \leq p,\ j=1,2, \cdots,r,
\end{eqnarray*}
where $\epsilon_{ij}$'s are iid $N(0,\sigma^2)$. We take $\sigma^2=1$ for
simplicity. We do not know which model is true. Also $p$ is large.

Assume now the unknown parameters $\beta_1,\beta_2,\cdots,\beta_q$ are
iid $N(0,c/r)$. Following Mukhopadhyay and Ghosh (2003) and Berger and
Pericchi (2001) we use r=1. This is now a new model with only two unknown
parameters $q$ and $c$, in place of the original $\beta_1,\beta_2,\cdots,\beta_p$.
The new PEB formulation is a considerable simplification compared with the
original version with $p$ parameters, $p$ large. The PEB approach is to estimate
$c$ from data and put a prior $\pi(q)$ on $q, \ 1 \leq q \leq p$. A typical
choice of $\pi(q)$ would be the uniform.

The likelihood of $(q,c)$'s, with $\beta_j$'s integrated out is
\begin{eqnarray*}
L(q,c)\propto (1+c)^{-q/2}\exp{\{-\frac{SS_q}{2(1+c)}-\frac{Y'Y-SS_q}{2}\}},
\end{eqnarray*}
where $\mathbf{SSq=\sum^q_{i=1}Y_i^2}$. \textbf{Let}
$\mathbf{\Lambda(q,c)=2\log{L(q,c)}+Y'Y}$, \textbf{then}
\begin{eqnarray*}
\mathbf{\Lambda(q,c)=\frac{1+c}{c}SS_q-\log{(1+c)}q.}
\end{eqnarray*}
Since $c$ is not known, one choice is to fix $q$ and maximize the conditional likelihood
with respect to $c$. The conditional mle of $c$ is given by
\begin{eqnarray*}
\hat{c}_q=\max{\{\frac{SS_q}{q}-1,0\}}.
\end{eqnarray*}

We now take $\pi(q)$ uniform on $1\leq q \leq p$. Then the PEB rule will choose
$M_q$ if $q$ maximizes $\Lambda(q,c)$ after replacing $c$ by $\hat{c}_q$, leading
to model selection criterion
\begin{eqnarray*}
\Lambda(q,\hat{c}_q)=SS_q-q(1+\log_{+}(\frac{SS_q}{q})),
\end{eqnarray*}
where $\log_{+}(\ )=\max{\{\log(\ ),0\}}$.

First fix $0<\alpha<1$ and define an $\alpha$-quantile model as defined by $M_q$
where $\pi(i+1\leq q|Y)\leq \alpha < \pi(i\leq q|Y)$.
The Bayes rule for best prediction chooses the model with smallest
dimension if $rc\leq 1$, and the $(c-1)/(2c)$ quantile model if $c>1$.

Since $c$ is not known we replace $c$ by $\hat{c}$ as discussed before and
plug in $\hat{c}$ for $c$. It can be shown this model is asymptotically as
good as an oracle for posterior prediction loss.
\textit{We assume we  use the least squares estimator
under each $M_q$. We get a posterior median rule if the predictor is the
Bayes estimator $E(\beta_i|q,Y)$ as in Barbieri and Berger (2000) later
published as Barbieri and Berger (2004).} The form of the predictor
(least squares or Bayes) has a substantial effect on model selection.
The optimality of AIC, discussed earlier, assumes the predictors are based
on least squares estimators of regression coefficients.

In their original paper George and Foster (2000), consider models,
\begin{eqnarray*}
M_{\gamma}\ : \ Y=X_{\gamma}\beta_{\gamma}+\epsilon,\ \epsilon
\stackrel{\mathrm{iid}}{\sim} N(0,\sigma^2),
\end{eqnarray*}
where $\gamma$ indexes all the subsets of $\ x_1,\cdots,x_p$, where $X_{\gamma}$ is the $n\times q_{\gamma}$ matrix whose columns are the $q_{\gamma}$ variables in the $\gamma$-th subset.

For model $M_{\gamma}$, the prior $\pi(\beta_{\gamma},\gamma)=\pi(\beta_{\gamma})\pi(\gamma)$. Here $\pi(\beta_{\gamma})$ is the Zellner-Siow prior for nonzero $\beta$'s penalizing the complexity  of the model and $\pi(\gamma)=w^{q_{\gamma}}(1-w)^{p-q_{\gamma}}$ is a binomial prior with $0\leq w \leq 1$, penalizing multiplicity. This may be the first time both complexity and multiplicity are recognized and penalized for model selection. George and Foster (2000) show how various model selection rules, including both AIC and BIC, can be derived in a unified way through this approach. However, this would ignore the goal of AIC as choosing a model $M_q$ that leads to best prediction.

\section{Consistency in High Dimensional Model Selection}\lel{s3}

Technically, consistency of model selection in high dimensional problems is usually very hard to 
prove, because our knowledge of the asymptotic behavior of the Bayes factor under a fixed $\theta$ 
is still quite meagre. For fixed (or slowly increasing) dimensions we have approximations like the 
BIC. These are no longer generally available in the high dimensional case.

In this section we try to do three things. We first discuss in detail a high dimensional model 
selection problem proposed by Stone (1978) and discussed in Berger et. al. (2003). Then we explore 
tentatively a weaker definition that may be both \textbf{more} satisfactory and often verifiable. The new definition 
reduces model selection to a test of two simple hypotheses. In at least one problem we show, we can 
get quite definitive results for a whole class of priors (Theorem 1).

We have also made a couple of general conjectures involving positive results on consistency 
(in the usual sense) under $M_1$. Consistency under $M_1$ is very important from the point of 
view of parsimony.

In Subsection 3.1 we try to initiate a change in our perspectives. We first attempt to pose 
a few general questions, make conjectures and suggest how they may be answered. It is no more 
than a modest attempt to inspire others to take up the challenges of consistency in model 
selection and, through that effort, get a better insight about high dimensional model selection.

\subsection{Generalized BIC for Stone's High Dimensional Example and Posterior Consistency}\lel{s31}

Common inference procedures for model selection in high dimensional examples
may behave in a very different way from those in low dimensional examples.
This subsection gives such an example, essentially due to Stone, studied in
Berger et. al. (2003).  This subsection is based on Berger et. al. (2003).
We try to initiate new work, as indicated earlier.

The following is a slight modification of a high dimensional linear model due
to Stone, showing (as in Stone's original example) that in high dimensional
problems, BIC need not be consistent but AIC is. \textbf{Let us consider}
\begin{eqnarray}\label{2ways}
y_{ij}=\mu_i+\epsilon_{ij},\ i=1,\ldots,p, \ j=1,\ldots,r, \ and
\ \epsilon_{ij} \ are \ i.i.d. \ N(0,1).
\end{eqnarray}
The two competing models are $M_1: \mu=0$ and $M_2: \mu \in R^p$. \textbf{The dimension $p\to \infty$.}

Berger et. al. (2003) show that the BIC is a very bad approximation to the
marginal under one of the models, and hence its inconsistency is unrelated
to consistency of Bayes factors in general. We reproduce a table (Table 1)
showing the poor approximation obtained from BIC and the much better
approximations obtained from the generalization, called GBIC, in Berger et. al.
(2003). A naive application of the relatively standard Laplace approximation due
to Kass, Wasserman and Pauler, denoted in Berger et. al. (2003) as $Lap_{KWP}$
provides a good approximation but not as good as the best. The best approximation
is provided by the rigorous Laplace approximation for the BF based on high
dimensional Cauchy,  given in the last column. But even this approximation
is poor near small values of $c_p$, i.e., when the alternative model ($M_2$)
is likely to be close to the null model (i.e., $M_1$).

\begin{table}[htbp]
	\centering
		\begin{tabular}{|c|c|c|c|c|c|}\hline
		$c_p$ & $\log{BF^c}$ & BIC & GBIC & $Lap^c_{KWP}$ & $\log{BF^{HD_c}}$\\\hline
.1 & -8.53 & -110.12 & -1.95 & 15.12 & -8.57\\
.5 & -3.82 & -90.12 & -1.95 & -2.26 & -3.908\\
1.0 & 6.03 & -65.12 & 5.71 & 5.55 & 5.92 \\
1.5 & 20.82 & -40.12 & 20.57 & 20.38 & 20.75\\
2.0 & 38.48 & -15.12 & 38.38 & 39.13 & 38.44\\
2.1 & 42.23 & -10.12 & 42.16 & 41.89 & 42.19\\
10.0 & 397.36 & 384.87 & 398.15 & 397.29 & 397.36\\\hline 		
			
		\end{tabular}
		\caption{Log (Bayes factor) and its approximations under the Cauchy prior}
\end{table}

We now turn to consistency issues, our main interest in this section.
Berger et. al. (2003) consider a family of priors including the popular
Zellner-Siow multivariate Cauchy prior and the Smooth Cauchy prior due to
Berger and Pericchi (1996). It turns out that the Bayes factor for the
multivariate Cauchy prior is consistent under both $M_1$ and $M_2$, but
is not consistent under the smooth prior. At the time this paper was written,
the paper pioneered a study of consistency of Bayesian model selection as well
as use of BIC to approximate a Bayes factor in high-dimensional problem.
In retrospect one understands consistency much better and some more insight as
 well as new results can be provided. We do this below.

We first note that the two models are nested, the parameter space under $M_1$ is
a singleton, namely it has only the point with all coordinates equal to zero,
and finally though this point is not contained in $\Theta_2$, it lies in
$\bar{\Theta}_2$. This is basically like testing a sharp null and with a
disjoint but not topologically separated \textbf{alternative}, i.e., we can not
find two disjoint open sets, one containing $\theta_0=0$ and the other
$\Theta_2$. Thus any open set containing $\theta_0$ will have non-empty overlap w
ith $\Theta_2$. For a discussion of the same point in Bayesian nonparametrics, see
Tokdar et. al. (2010).

That, the posterior for the Zellner-Siow prior (denoted as Z-S prior or $\pi_c$) is consistent inspite of lack of separation of $M_1$ and $M_2$ may be intuitively explained by examining the structure of this prior as well as that of the Smooth Cauchy prior (denoted as $\pi_{sc}$). The Z-S prior may be represented as a mixture of multivariate normals with a gamma prior for the common precision parameter $t=\frac{1}{\sigma}$,
eqn. (6) of Berger et. al.
\begin{eqnarray*}
\lefteqn{\pi_{c}(\mu)= \frac{\Gamma((p+1)/2)}{\pi^{(p+1)/2}}(1+\mu'\mu)^{-(p+1)/2}}\\\nonumber
&&=\int_0^{\infty}\frac{t^{p/2}}{(2\pi)^{p/2}}e^{-(t/2)\mu'\mu}\frac{1}{\sqrt{2\pi}}e^{-t/2}t^{-\frac{1}{2}}dt.
\end{eqnarray*}
The gamma mixing measure puts positive mass near precision parameter $t=\infty$, i.e., in the $\sigma$-space, near $\sigma=0$, making even small differences detectable (under both $M_1$ and $M_2$). Hence $\pi_c(\mu)$ is a consistent prior. On the other hand
\begin{eqnarray*}
\pi_{sc}(\mu)=\int_0^1\frac{t^{p/2}}{(2\pi)^{p/2}}e^{-(t/2)\mu'\mu}\frac{1}{\pi\sqrt{t(1-t)}}dt
\end{eqnarray*}
has inconsistent posterior (under $M_2$) because it is supported on the set (0,1),
with $t=\infty$, i.e. $\sigma=0$ (in the $\sigma$-space) not in the support of the
mixing measure.
Some comments on these facts and proofs in Berger et. al. (2003) are also in order.

Consider the general family of priors
\begin{eqnarray}\label{priors}
\pi_g(\mu)=\int^{\infty}_{0}(\frac{t^{p/2}}{(2\pi)^{p/2}}e^{-(t/2)\mu'\mu})g(t)dt.
\end{eqnarray}
All priors with the general structure given by Berger et. al. (2003) in \textbf{(4)}
ensure posterior consistency under $M_1$ even if $g(t)$ is strictly positive
only on $(0,T]$, for arbitrary $T>0$. The condition $g(t)>0$ on $(0,\infty)$
is needed for consistency under $M_2$, but not under $M_1$. (In the proof of
Theorem 3.1 in Berger et. al (2003), without this assumption the set $S_{\epsilon}$ is
empty but that does not matter under $M_1$.) Thus the smooth Cauchy with support of
$g$ equal to $(0,1]$ has consistent posterior under $M_1$, and under $M_2$ it is
inconsistent if $0<\tau^2<2\log(2)-1$, where $\tau^2=\lim_{p\to \infty}
\frac{1}{p}\sum{\mu_i^2}$.

A similar result holds if the support of $g$ is $(0,T]$. The most
interesting fact that emerges is the rather general consistency theorem
under $M_1$. For those of us who believe in parsimony, this is a very pleasant fact.

In Scott and Berger (2010), to be discussed in Subsection 5.4, we note a similar
fact under the global null, provided the global null is given a positive prior
probability. One gets posterior consistency under this model by straightforward
applications of Doob's Theorem (Ghosh and Ramamoorthy, 2006, p22), as in
the proof of similar results for Bayesian Nonparametric models by Dass and Lee (2004).

\textbf{Conjecture}. All this suggests a general consistency theorem under the
global null for linear models remains to be discovered and that one would need
to try to prove a general version of Doob's theorem for independent but not
identically distributed r.v.'s. A simple example appears in Choi and Ramamoorthi
(2008).

So far in this paper we have been using consistency as defined in Section 2.1.
We now define a new notion of consistency, which is simply consistency of Bayes
factors under the two marginals of $\{X_i\}$ under $M_1$ and $M_2$.

\textbf{Definition of $\mathbf{(P_1,P_2)}$-consistency} Let $P_1$ and $P_2$ be the infinite-dimensional marginal distributions of $X_i$'s under $M_1$ and $M_2$. We say consistency holds, iff under the true infinite-dimensional marginal distributions $P_1,\ P_2$, $\log{BF_{21}}\to -\infty$ and $\infty$ in probability. We call this $(P_1,P_2)$-consistency to distinguish from usual consistency.

This should be particularly attractive, if the data $(X_1,\cdots, X_n)$
have been generated as assumed in Scott and Berger (2010), i.e., under $P_1$
or $P_2$ as true, not an unknown $\theta$. If this model is correct, a large
enough data is expected to choose the correct model with very high probability.

We prove $(P_1,P_2)$-consistency for Stone's example.

\begin{theorem}
Assume the general family of priors in equation (\ref{priors})
with $0<t<\infty$ w.p. $1$ under $g(t)$. Then $(P_1,P_2)$-consistency holds.
\end{theorem}

\textbf{Proof}
Conditionally, for fixed value of the precision parameter $t$, under
$M_2$, $\frac{1}{p}\sum^p_{1}\bar{X}^2_i\to \frac{1}{r}+\frac{1}{t}$ a.s.
Since $t>0$ with probability one, under $P_2$ (obtained by integrating out $t$),
\begin{eqnarray*}
\lim \frac{1}{p} \sum^p_{1}\bar{X}^2_i >\frac{1}{r}, \, \, a.s.
\end{eqnarray*}
On the other hand, under $P_1$,
\begin{eqnarray*}
\lim \frac{1}{p} \sum^p_{1}\bar{X}^2_i \to \frac{1}{r}, \, \, a.s.
\end{eqnarray*}

The above facts show $P_1$ and $P_2$ are orthogonal (i.e., $P_1$ and $P_2$ are supported on disjoint subsets of the sample space $\{X_i\}$) and hence, by standard facts about likelihood ratios (Kraft, (1955)) the Bayes factor
\begin{eqnarray*}
\log{BF_{21}}&\to& \infty \ \ \ a.s. (P_2), \qquad
\log{BF_{21}}\to -\infty \,\,\, a.s. (P_1),
\end{eqnarray*}
proving posterior consistency under both $M_1$ and $M_2$. $\Box$

It seems plausible to us that just as consistency under $M_1$ may hold rather 
generally (by Doob's theorem) for linear models, posterior consistency under $P_2$ 
may also be true rather generally. Moreover even when the question cannot be settled 
theoretically, a Bayesian simulation as in Scott and Berger (2010) will throw light 
on whether consistency is to be expected or not. All this will bypass the need to have 
good Laplace approximations to $m_1(\mathbf{x})$, $m_2(\mathbf{x})$ in high dimensional cases.

The general linear model, of which Berger et. al. (2003) is a very simple special case, 
has been studied in a greater detail in Liang et al. (2008) from the point of view of choice of 
new priors, calculation of Bayes factor and consistency (for the fixed $p$ case). However, posterior 
consistency for the high dimensional case doesn't seem to be have been studied.

On the other hand Moreno et. al. (2010) study posterior consistency for some general 
high dimensional regression problems when intrinsic priors are used. The results are 
quite interesting but the formulations of consistency are somewhat different.

\section{Cross-Validatory Bayes factor}\lel{s4}

\subsection{General Issues}\lel{s41}
It has been known for quite some time that 
Bayesian estimation of parameters or prediction of future observations is quite robust with respect
to the choice of prior while model selection (for 0-1 loss) based on Bayes factors is not robust.
Draper and Krnjajic (2010) discuss instability of Bayes factors and suggest replacing
them with cross-validatory Bayes factors.

Very roughly speaking, in estimation one uses diffuse improper or diffuse proper
priors, so that most of the information in the posterior come from the data, not
the prior. In particular the undetermined constants in an improper prior gets
canceled because it appears in both the numerator and denominator of the basic
Bayes Formula:
\begin{eqnarray*}
p(\theta|x)=\frac{c.p(\theta)p(x|\theta)}{\int{c.p(\theta)p(x|\theta)d\theta}}\, .
\end{eqnarray*}

On the other hand testing procedures or model selection methods based on Bayes \textbf{factor} do
not have these good properties, see e.g., Ghosh and Samanta (2002) and 
Ghosh et. al. (2006). A standard way of solving both problems is to use data
based priors as follows. One uses a part of the data, say a vector $x_{k}$ to
calculate the posterior for each model which is then used as a prior for the
corresponding model. For each model the corresponding data based prior is then
combined with the remaining data set $x_{-k}$ consisting of the 
remaining $(n-k)$ $x_i$'s, to produce a marginal based on $x_{-k}$. With
these modifications, both difficulties, namely the appearence of the arbitrary
constant in an improper prior and lack of robustness with respect to the prior
disappear. All of the new priors are really posteriors and, hence, robust if based
on substantial data.

There is a huge literature on these cross-validatory Bayes factors. We mention a few,
based partly on their importance and partly on our familiarity, Geisser (1975),
Berger and Pericchi (1996), Ghosh and Samanta (2002), O'Hagan (1995), Chakrabarti
and Ghosh (2007) and Draper and Krnjajic (2010). Berger and Pericchi (1996, 2004)
adopt the same procedure but condition w.r.t., what they call the minimal training
sample, that is a smallest subsample such that conditioning w.r.t. it makes the
posterior proper. Then averages are taken over minimal training samples. We note
in passing, that this method has actually led to construction of objective priors,
 which Berger and Pericchi (2004) \textbf{call} intrinsic priors. The papers of Berger and
Pericchi (2004) and \textbf{other} colleagues \textbf{provide} many details and interpretations. \textbf{We} need to know only the basic facts stated above.
\textit{Berger and Pericchi are not trying to get a stable Bayes factor,
they are trying to get a \textbf{Bayes} factor which may be called ``objective''.}

These new Bayes factors are more stable but raise a new issue which is still not
fully understood. How much of the data should be used to make the prior stable
(by computing the posterior and treating the posterior as a data based prior)
and how much of the data should be used for inference? We first heard this
question from Prof. L. Pericchi. This is a deep and difficult question.
It has been discussed in Chakrabarti and Ghosh (2007). In the next subsection
we turn to the problem involving what level of cross-validation should be chosen
if we are in the M-closed case.

The cross-validatory Bayes factors have a long history, suggested by
Bernardo and Smith (1994). They were extended by Gelfand and Dey (1994),
who in turn had drawn on Geisser (1975). Consistency issues were studied
by Chakrabarti and Ghosh (2006), assuming these Bayes factors as given -
a major motivation was to throw some light on Pericchi's question but not
settle any of the other basic normative issues, specially in the context
of selecting a true model or one close to it in some sense. Many interesting
alternative approaches were suggested by discussants of Chakrabarti and Ghosh (2007),
namely Lauritzen, Pericchi, Draper, Vehtari and others, which are still not explored.
In the next subsection we content ourselves with revisiting a partly heuristic
treatment of the high dimensional M-closed nested case, i.e. , with either $M_1$
or $M_2$ true. We provide a relatively simple proof of consistency under $M_2$
and $M_1$. A partly heuristic proof of consistency under $M_1$ is given in
Chakrabarti and Ghosh (2007). Generalizing this result we also suggest how the
sample size $r$ is to be allocated between stabilizing the posterior and model
selection under more general assumptions than Chakrabarti and Ghosh (2007).

A very recent paper is Draper and Krnjajic (2010), who suggest the
cross-validatory Bayes factors be used for model choice in the M-closed
case to ensure stable inference. Though they report the results of a few
simulations which are promising, one would need a much more extensive study of
complex varying dimensions and different sizes before drawing firm conclusions.

We make a few tentative comments about their work. Draper and Krnjajic (2010) seem to be 
giving up the usual Bayes factors or replace them with cross-validatory Bayes factors. 
Also they seem to make simplistic assumptions about asymptotics in model selection. 
When $n$ goes to infinity, $p$ will usually tend to infinity, but not necessarily so. 
Moreover the rate of growth of $p/n$ can vary a lot, leading to very different kinds of 
asymptotics. In particular if $p$ tends to infinity sufficiently slowly, the results will be 
like those for fixed $p$ and $n\to \infty$.

It will be interesting to compare the cross-validatory BF of Draper and Krnjajic (2010) and 
Intrinsic Bayes factor in the same problem. Also, our general view for fixed dimensional parameter 
and moderate n, is that the cross-validatory BF does not differ much from the usual BF. 
The following heuristics might clarify why this is likely. We consider 
$X_1,\ldots, X_n\sim N(\theta,1)$ and competing models $M_1: \theta=0$ and $M_2:\theta>0$, 
with a standard prior $N(\mu,\sigma^2)$ for $\theta$. The leading term of the marginal under 
$M_1$, by Laplace approximation can be written as,
\begin{eqnarray*}
-\frac{1}{2}\sum_{j=1}^n {(x_j-\bar{x})^2}.
\end{eqnarray*}
A heuristic cross validatory replacement of this will be 
\begin{eqnarray*}
-\frac{1}{2}\sum_{j=1}^n {(x_j-\bar{x}_{-j})^2}, \qquad \bar{x}_{-j}=\sum_{i\ne j} x_i/(n-1).
\end{eqnarray*}
 The fact $\sum{(x_j-\bar{x}_{-j})^2}=\sum{(x_j-\bar{x})^2}(\frac{n}{n-1})^2$ has the effect 
of reducing the marginal under alternative, making it conservative under $\theta=0$. 
The leading term of cross-validatory BF increases a bit, suggesting it will work better. 
This has been confirmed by using cross-validatory BF and BF, with simulation studies. 
We have seen that under $M_2$ they differ by a very small value in log-scale, but that 
the small difference plays a critical role under $M_1$. Under $M_1$, $CVBF_{12}$ is greater 
than one for 80\% cases choosing the correct model compared to none of them in case of $BF_{12}$, 
but we are looking at cases where the Bayes factor is very close to one.

We can justify the above claims for the cross-validatory Bayes factor of
Draper and Krnjajic (2010), who replace the marginals in the Bayes factor by
\begin{eqnarray}\label{eq1}
\frac{1}{n}\sum^n_{i=1}\log{p(y_i|y_{-i},M_j)}.
\end{eqnarray}
Justification follows from a straight forward application of results in
Mukhopadhyay et. al. (2005). Note that the above CVBF is the Bayes factor
defined in the last reference with $s=n-1$. Then the identity in
Mukhopadhyay et. al. (2005, eq. 11) reduces equation (\ref{eq1}) to
\begin{eqnarray*}
-\frac{1}{2}\sum^n_{i=1}{(x_i-\bar{x})^2}-\frac{n-1}{2}\log{2\pi} -\frac{1}{2}\log{n}
\end{eqnarray*}
The proof of the equation 11 in Mukhopadhyay et. al. (2005) is given in the
appendix of that paper.

\subsection{How to choose the size of cross-validation : some preliminary results}\lel{s42}

Following Chakrabarti and Ghosh (2007), we will review the high-dimensional
normal linear model setup as described in
equation (\ref{2ways}), with the same \textbf{competetive} models $M_1$ and $M_2$ described
there.
The study in Chakrabarti and Ghosh (2007) was done under the Zellner-Siow prior,
but our results hold for any other priors under the following assumption.
We assume that\\[.15cm]
{\bf Assumption 1.} \,\, $ \pi(\hat{\mu}_k)-\pi(\hat{\mu}_r)=o_p(1),$ as
$ k\to \infty,$ $ r\to \infty$,
where $\mathbf{\pi(\hat{\mu}_k)}$ and $\mathbf{\pi(\hat{\mu}_r)}$ are the prior density evaluated at the mle of
$\mu$ depending on $k$ and $r$ replicates.
\vskip .15cm 
For cross-validatory Bayes factor we use $k$
out of $r$ replicates for each $\mu_i$  to make the prior proper. A formal
definition appears in equation (\ref{19}) below. Here we will try to prove the
consistency of
a proxy to $CVBF_{21}$ obtained by using the popular KWP-Laplace approximation for
high dimensional problems. Under Assumption 1, below we apply this approximation
to both $BF^r_{21}$ and $BF^k_{21}$ with $p$ as dimension of parameter space and
$r$ and $k$ as sample size.
\begin{eqnarray}\label{19}
\log{CVBF_{21}}&=&\log{\frac{BF^r_{21}}{BF^k_{21}}}
=\frac{p}{2}[(rC_p-kC'_p)-\log{\frac{r}{k}}]+o_P(1)\\
&=&\log{CVBF^{ps}_{21}}+o_P(1), \notag
\end{eqnarray}
where $C_p=\frac{1}{p}\sum^p_{i=1}[\frac{1}{r}\sum^r_{j=1}{y_{ij}}]^2$
and $C'_p=\frac{1}{p}\sum^p_{i=1}[\frac{1}{k}\sum^k_{j=1}{y_{ij}}]^2$.
Here ``ps" stands for pseudo.

We will prove posterior consistency under both models $M_1$ and $M_2$ for
$CVBF^{ps}_{21}$. To prove the consistency under model $M_2$, we assume as
in Chakrabarti and Ghosh (2007), the following.\\[.1cm]
{\bf Assumption 2.} $\lim_{p\to \infty} \frac{1}{p} \sum^p_{i=1}\mu^2_i=\tau^2>0$,
under $M_2$.

\begin{theorem} Let $p,\ k, \ r \to \infty$. Under Assumptions 1 and 2,
for both $M_1$ and $M_2$, $\log{CVBF^{ps}_{21}}$ chooses the correct model, 
with probability tending to one, for all $0 \leq c<1,$ when  $k/r\to c $.
\end{theorem}
\textbf{Proof}
Suppose first, model $M_2$ is true, then for fixed value of p,
\begin{eqnarray*}
\log{CVBF^{ps}_{21}}&=&\frac{p}{2}[(rC_p-kC'_p)+\log{\frac{k}{r}}]\\
&=&\frac{p}{2}[(rC_p-kC'_p)-\log{\frac{rC_p}{kC'_p}}+\log{\frac{C_p}{C'_p}}]
\end{eqnarray*}
The function $f(x)=x-\log{x}$ is increasing for $x>1$. Using Assumption 2, 
when $k\to \infty, 
r\to \infty$, we have for sufficiently 
large $p,\, r$ and $k$, $P(rC_p>1)\to 1$, $P(kC'_p>1)\to 1$. Then $\frac{rC_p}{kC'_p}>1$ 
implies $(rC_p-kC'_p)-\log{\frac{rC_p}{kC'_p}}>0$. By definition of $C_p$ and $C'_p$, and also 
using Assumption A2, under $k\to \infty$ and $r \to \infty$
\begin{eqnarray*}
\lim_{k\to \infty, r\to \infty} \frac{p}{2}\log{\frac{C_p}{C'_p}}=\frac{p}{2}o_P(1)
\end{eqnarray*}
which implies $\frac{C_p}{C'_p}=1+o_P(1)$. Hence, we obtain 
\begin{eqnarray*}
\lim_{k\to \infty, r\to \infty}\frac{rC_p}{kC'_p}=\frac{1}{c}.
\end{eqnarray*}
Now here $\frac{1}{c}>1$ for any $0 \leq c<1$, completing the proof of consistency under $M_2$.

We know,
\begin{eqnarray*}
\lefteqn{(rC_p-kC'_p)}\\
&&=(\frac{k^2}{r}-k)\frac{1}{p}\sum{\bar{y}^{'2}_i}+r(1-\frac{k}{r})^2\frac{1}{p}\sum{\bar{y}^{''2}_i}+2k(1-\frac{k}{r})\frac{1}{p}\sum{\bar{y}'_i\bar{y}''_i}\\
&&=(\frac{k}{r}-1)\frac{1}{p}\sum{(\sqrt{k}\bar{y}'_i)^{2}}
-(\frac{k}{r}-1)\frac{1}{p}\sum{(\sqrt{r-k}\bar{y}''_i)^{2}}\\
&& \hskip 1in +2\sqrt{\frac{k}{r}(1-\frac{k}{r})}\frac{1}{p}\sum{\sqrt{k}\bar{y}'_i\sqrt{r-k}\bar{y}''_i}
\end{eqnarray*}
Now assuming model $M_1$ is true and $r\to \infty$, $r-k\to \infty$, $\frac{k}{r}\to c$, 
for fixed value of $p$, we get,
\begin{eqnarray*}
\lefteqn{\lim_{r\to \infty, k\to \infty}rC_p-kC'_p}\\
&&=(c-1)[\frac{1}{p}\sum{W_j}-\frac{1}{p}\sum{U_j}]
+2\sqrt{c(1-c)}\frac{1}{p}\sum{w_i u_i},
\end{eqnarray*}
where $W_j \sim \chi^2_1$, $U_j\sim \chi^2_1$, $w_i\sim N(0,1)$ and $u_i \sim N(0,1)$ for $\forall\ 1<j<p$. Hence,
\begin{eqnarray*}
\lefteqn{ \lim_{k\to \infty, r\to \infty}\frac{p}{2}[(rC_p-kC'_p)-\log{\frac{r}{k}}]}\\
&&=\frac{p}{2}[(c-1)[\frac{1}{p}\sum{W_j}-\frac{1}{p}\sum{U_j}]+2\sqrt{c(1-c)}\frac{1}{p}\sum{w_i u_i}+\log{c}]\\
&&=\frac{p}{2}[o_p(1)+\log{c}]\sim \frac{p}{2}\log{c}<0,
\end{eqnarray*}
for $0\leq c<1$, showing the consistency of $\log{CVBF^{ps}_{21}}$ under $M_1$. $\Box$

\textbf{Remark} Extending the results proved in Chakrabarti and Ghosh (2007), we have shown 
the consistency of the cross-validatory Bayes factor under both $\mathbf{M_1}$ and $\mathbf{M_2}$ for 
$0\leq c <1$, when $\frac{k}{r}\to c$. Our proof is valid for any prior distribution satisfying 
the Assumption 1, which includes the Zellner-Siow prior used in Chakrabarti and Ghosh (2007) 
and many other commonly used. The proof also shows the smaller the value of $c$ (i.e. smaller the $k$) 
we get a larger $CVBF_{21}$ under $\mathbf{M_2}$ and a smaller $CVBF_{21}$ under $\mathbf{M_1}$. This suggests one 
should have a relatively small value of $c$, $c=0$ would be the best. This supports the choice of 
minimal training sample as in Berger and Pericchi (2004) but seems to contradict the conjecture of 
Chakrabarti and Ghosh (2007) that in high dimensional case c should tend to a positive constant. 
A possible explanation is that $r\to \infty$, $k\to \infty$ but $p$ tends to infinity at a slower 
rate than $k$ and $r-k$, this prevents it from becoming a real high-dimensional problem. 
Then further study is needed.

\section{Multiple Testing : General Issues}\lel{s5}

In recent years multiple testing  has emerged as a very important problem in 
statistical inference, because of its applicability in understanding large data sets. 
One of the major fields of applications is bioinformatics, where multiple testing is 
extensively applied for the analysis of gene expression, proteomics or genome wide 
association studies (GWAS).
The following subsections deal with Full Bayes, Empirical Bayes (EB) and the 
classical approach to multiple testing based on FDR (False Discovery Rate). 
Multiple testing can be reformulated as model selection. This has been known in the 
early literature on this subject, see for example Hodges (1956).

\subsection{The Full and Empirical Bayes Approaches for multiple testing}\lel{s51}

The Empirical Bayes approach is very popular in the context of analyzing the microarray 
data. It is well understood that due to a very small number of replicates standard 
maximum likelihood estimates of the variance of the individual gene expression are 
very imprecise.  Therefore the Bayes hierarchical model is often used and the standard 
deviation for each of the respective t-statistics is estimated based on the Empirical 
Bayes approach. The method is implemented in the package LIMMA (Smyth (2004)), which 
became a standard tool for practical microarray analysis.

 More advanced Empirical Bayes approach is proposed  in Datta and Datta (2005), 
who use the kernel density estimate to estimate the marginal density of 
$z_i=\Phi^{-1}(p_i)$, where $\Phi(\cdot)$ is the cdf of the standard normal 
distribution and $p_1,\ldots,p_m$ are p-values for consecutive test statistics. 
Then, assuming that $z_i\sim N(\theta_i,1)$, where $\theta_i=0$ corresponds to 
the null hypothesis, the authors calculate the Empirical Bayes estimate for 
$\theta_i$ and use resampling to decide on the corresponding threshold. 
As shown in Datta and Datta (2005), in many cases their method offers a 
substantially larger power than the popular Benjamini-Hochberg procedure 
(BH, Benjamini-Hochberg (1995)), which however happens at the price of some 
increase of FDR.

Many contemporary multiple testing procedures are  based on a two component mixture 
model, used e.g. in Berry (1988),  Efron et al. (2001) or Efron and Tibshirani (2002). 
Such a mixture model assumes that test statistics
$X_1,\ldots, X_m$ are iid rv's and their marginal cdf $F(x)$ can be modeled as
\begin{equation}\label{model}
F(x)=(1-p)F_0(x) + p F_A(x)\;\;,
\end{equation}
where $F_0(x)$ and $F_A(x)$ denote the cdfs of the null and alternative distributions, 
respectively, and $p$ is the expected proportion of alternatives among all tests.

A variety of Empirical Bayes methods, both parametric and nonparametric, have been 
proposed for the estimation of the unknown functions and parameters in
(\ref{model}) (see, e.g., Efron and Tibshirani (2002), Johnstone and Silverman (2004), 
Storey (2007), Jin and Cai (2007), Bogdan et al. (2007, 2008a), Efron (2008) or Cai 
and Jin (2010)).  These estimates are subsequently used to approximate multiple 
testing rules controlling False Discovery Rate (FDR) or local FDR  (Efron and 
Tibshirani (2002)),  minimizing the estimation error (Johnstone and Silverman (2004)) 
or  for approximation of classical oracles, aimed at maximizing the power while 
controlling some measures of the type I error or FDR (see, e.g., Storey (2007) or 
Sun and Cai (2007)).
In Bogdan et al. (2007, 2008a) Empirical Bayes estimates are used to 
approximate Bayesian oracle, aimed at minimizing the Bayes risk.

It is also interesting to observe that the Benjamini-Hochberg procedure (BH,Benjamini and Hochberg, 1995), aimed at controlling FDR, shares some similarities with the Empirical Bayes methods. Namely, BH  can be understood as the plug-in rule to control the, so called, Bayesian False Discovery Rate (e.g. see Efron and Tibshirani (2002) or Bogdan et al. (2011)). In BH the cdf of the mixture is estimated by the empirical cdf and in this way BH borrows strength from information provided by all data points, similarly as it is done in Empirical Bayes methods. Furthermore, in Bogdan et al. (2011) and Frommlet et al. (2011) it is shown that under  sparsity ($p\rightarrow 0$) BH  has some asymptotical optimality properties in the context of minimizing the Bayes risk (for more details see Section 5.2), even though it is a non-Bayesian mutiple test.

The full Bayes approach to multiple testing is discussed e.g. in Scott and Berger (2006), who consider the parametric setting under which both the conditional  distribution of the test statistic as well as the prior for the effect size are normal. The fully Bayes nonparametric multiple testing procedure based on MCMC algorithm of Escobar and West (1995) is investigated e.g. in Bogdan et al. (2008a).  In a general article Scott and Berger (2010) explain  how the fully Bayes procedures for model selection correct for multiplicity by using a fixed prior distribution on $p$.

\subsection{Sparsity and limitations of multiple testing procedures}

In many applications of multiple testing it is assumed that the proportion of alternatives among all tests $p$ is very small. In the asymptotic context this assumption is often summarized by letting $p\rightarrow 0$ as the number of tests $m$ goes to infinity. Recently Abramovich et al. (2006), Bogdan et al. (2011) and Frommlet et al. (2011b) analyzed the asymptotic properties of the multiple testing procedures under the sparsity. Bogdan et al. (2011) and Frommlet et al. (2011b) have shown that the limiting power of the Bayes oracle is larger than zero only if a signal magnitude is large enough. In Bogdan et al. (2011) the signal magnitude $u$ is measured by the ratio of variances of the test statistics under the alternative and null distribution. It is proved  that in case when the number  $n$  of replicates used to calculate each of the test statistics remains constant as $m\rightarrow \infty$ then the signals on the verge of detectability satisfy $u \propto \log(\delta/p)$, where $\delta$ is the ratio of losses for type I and type II erros. In case when $p\propto m^{-\beta}$, for some $\beta>0$, and $\log \delta=o(\log m)$ this condition can be simplified to $u \propto \log m$.  In Frommlet et al. (2011b) it is further shown that  when $u$ is fixed and $n\rightarrow \infty$  then signals are on the verge of detectability if $n \propto \log m$.
Bogdan et al. (2011) and Frommlet et al. (2011b) also prove that for $\beta \in (0,1]$, signals on the verge of detectability and some mild conditions on the ratio of losses, the ratio of the Bayes risks of the optimal Bayes oracle and the Benjamini-Hochberg procedure at a fixed FDR level $\alpha \in (0,1)$ converges to 1.  This result is weaker than the sample space similarity but it illustrates the asymptotic optimality properties of BH with respect to the Bayes risk and complements findings of Abramovich et al. (2006), who prove that the hard thresholding rule based on BH has some asymptotic minimax properties with respect to the  risk of the estimation of the unknown vector of means.

\subsection{Frequentist properties}

Johnstone and Silverman (2004) analyze the  risk of estimation of the  vector of means by  different thresholding rules based on multiple testing procedures. They are mainly concentrating on the Empirical Bayes procedure, which assumes that the null and alternative distributions of the test statistics are known and uses the method of the maximum likelihood to estimate  $p$. They show that under the sparsity and some assumptions on the null and alternative distributions, their proposed EB thresholding rule has the same asymptotic minimax properties as BH. However, their simulation study shows  that EB procedure is more robust than BH for denser signals.
Similar conclusions can be drawn based on the simulation study of
Bogdan et al. (2008a), who compare different multiple testing procedures with respect to some other frequentist properties, like the misclassification probability, FDR or the power. They show that unless $p$ is very small BH is usually outperformed by Full or Empirical Bayes procedures. Bogdan et al. (2008a) also demonstrate serious  problems of the naive Empirical Bayes estimates in case when $p$ is very small and the parameters of the alternative distribution are not known.
It turns out that in this situation the Kullback-Leibler distance between some mixture densities with  very small values of $p$  and the densities  with $p$  close to 1 may be very small. This phenomenon is related to the problem of nonidentifiability of the mixture for $p=0$ and obviously results in a very bad performance of the maximum likelihood estimates and the corresponding multiple testing procedure. Similar problems were also discussed in George and Foster (2000), who use the Empirical Bayes approach to estimate the mixture parameters in the context of model selection. Bogdan et al. (2008a) solved this problem by using an informative beta prior on $p$, suggested in Scott and Berger (2006), and replacing the maximum likelihood estimate with the approximation to the posterior mode. The results reported in Bogdan et al (2008a) show that the proposed Empirical Bayes procedures performs well and its frequentist characteristics closely approximate the characteristics of the Fully Bayes procedure.

\subsection{Approximation of posterior probabilities}\lel{s54}

In Scott and Berger (2010) the Full and Empirical Bayes approaches 
are compared with respect to posterior model probabilities. This article 
is a major contribution which points out and clarifies a multitude of issues 
and advantages of handling them in a Bayesian rather than Empirical Bayesian 
way. Indeed, one of the major messages in the paper is that the two approaches 
may be very different, unlike the case of estimation and contrary to the usual 
view that the second is an approximation to the first.

Scott and Berger (2010) stress the fact that the standard type II mle of $p$ can 
attain 0 or 1 with a positive probability. In these cases the posterior 
probabilities of the total null or the full model, respectively, are equal to one. 
This obviously does not reflect properly the  uncertainty in the data and 
differs substantially from the answer provided by the fully Bayes 
solution.  Scott and Berger (2010) show that due to these problems the 
expected Kullback-Leibler divergence  between FB and EB posteriors for the variable selection  is infinite.

Most Bayesians who think of EB as a convenient approximation to full Bayes, would certainly 
be ready to delink Type 2 mle's from EB. However, we have not checked if the Bayesian 
penalty brings EB estimate of $p$ and selected variables closer to the full Bayes scenario, 
which may need further thought. These points are easily checked via simulation.

Another interesting point is the full Bayesian simulation of data, in which the
parameters change at each draw. It would be interesting to see if this tends to
heighten the difference between Full Bayes and EB. That the difference will not 
vanish if, as usual, the parameters are kept fixed during simulations, is clear
from the analysis based on the real data, but a better estimate of sparsity and
comparing only the inference given the ``actual" data may mitigate the difference.

We now summarize and explore what remains similar and what is not in
the two approaches. As explained in Bogdan et. al. (2008a, 2011), under additive losses, 
the EB multiple test approximates well the Bayes oracle, as usually defined for 
example in machine learning, see Wasserman (2003). This further suggests the 
Full Bayes multiple test and the EB multiple test are likely to agree. 
On the other hand very non-linear quantities like the marginal probability of 
a model $M_{j}$, on which some variables are included and some not, is 
likely to be very different in Full Bayes and EB.
Take for example the general integral (2) in Scott and Berger (2010) with
$p(M_{\gamma})$ given by (8) or (9) in their paper. In general naive EB will not provide
good approximations. Whether simple but good approximations exist is an
interesting question.

\subsection{Application to high dimensional multiple regression}

The relationship between multiple testing and the problem of model selection
in high-dimensional multiple regression  has been extensively discussed e.g. in
George and Foster (2000), Abramovich et al. (2006), Scott and Berger (2010) or
Frommlet et al. (2011b). Abramovich et al. (2006) and Frommlet et al. (2011b)
develop new model selection criteria by a direct analogy with FDR methods for
multiple testing. \.Zak-Szatkowska and Bogdan (2011) further explain the
relationship between the penalties of these criteria and priors for the
model dimension in the full Bayesian analysis, while Frommlet et al. (2011b)
prove their asymptotic optimality in the context of minimizing the
Bayes risk under orthogonal designs.  George and Foster (2000) propose
two new model selection criteria based on the Empirical Bayes estimates
of the mixture parameters and illustrate their good properties with the
simulation study. On the other hand, Scott and Berger (2010) point at the
advantages of the fully Bayes approach for high-dimensional multiple regression.

While in principle one can  easily use the Bayes approach and compare different models
by calculating their posterior probabilities, there still remains a challenging
problem of the search over a huge space of possible solutions.
In  case  when the number of regressors is so large that the standard optimization
or search procedures are not feasible, Fan and Lv (2008)  recommend to perform
an initial screening based on the marginal correlations between explanatory
variables and the response variable (sure independence screening, SIS).
Fan and Lv (2008) prove that under some conditions on the random design
matrix and the number of true regressors  SIS  will contain all true
regressors with probability converging to 1. They also present results
of some simulation studies illustrating good properties of this procedure.

However, a recent simulation study of Frommlet et al. (2011a) shows that in the
practical problems of genowe wide association studies (GWAS), SIS may lead to a
substantially wrong ranking of predictors. It is important to note that in case
of GWAS the rescaled random matrix of marker genotypes does not have a spherically
symmetric distribution, which is one of the assumptions of Fan and Lv (2008).
But, the discussion in Frommlet et al. (2011a) suggests that the similar problem
would appear also in case where predictors were the realizations of independent
normal variables, and  its real source  is the number of true regressors.
In simulations reported by Frommlet et al. (2011a) the number of true regressors
was equal to 40, which substantially exceeds the number of regressors simulated
by Fan and Lv (2008). In Frommlet et al. (2011a) the number of individuals $n$
was set to 649 and the total number of predictors was equal to 309788.
The values of regression coefficients were chosen in such a way that
the power of detection of true regressors by the model selection criterion mBIC2,
defined in Frommlet et al (2011b), almost uniformly covered $[0,1]$. It turns out
that for this choice of parameters the summary effect of small random sample
correlations between  the  gene of interest and  remaining  (39 or 40) causal
genes may create a large spourious correlation between this gene and the trait,
and this random component  can easily become a dominating part of the marginal
correlation. While this is true that with increasing $n$ the random correlations
would gradually disappear and the true genes could be appropriately ranked, one
should recognize that simulating 40 causal genes is just a toy approximation for
a true complexity of GWAS. It is currently believed that  complex, quantitative
traits are often influenced by a huge number of very weak genes, so called polygenes.
In this situation the basic assumption of Fan and Lv (2008) that the true number of
regressors is smaller than $n$ is no longer satisfied. In this case the task of
localizing important genes relies rather on identifying   large outliers than
testing that a given regression coefficient is zero. We believe that under this
scenario ``strong'' genes can be appropriately identified by SIS only if their
individual effects exceed the summary
effect of polygenes. Developing precise conditions for the identifiability of
such strong regressors remains an interesting topic for future research.

\subsection{Comparing more than two experimental conditions}

In many real life problems scientists observe certain characteristics over a  period of time or under a set of ordered conditions. Then the major point of interest is detecting some significant trends. In case when many characteristics are observed at the same time this leads to a multiple testing problem, where one needs to control both the number of falsely detected nonzero trends as well as the number of wrongly identified patterns. The statistical framework for these types of problem is discussed in Guo et al. (2010), who consider it as a generalization of the multiple testing problem of detecting both the effect and its sign, discussed in Benjamini and Yekutieli (2005). The popular measure to characterize the performance of the directional testing procedures is the mixed directional FDR (mdFDR), introduced in Benjamini et al. (1993). Benjamini and Yekutieli (2005) proposed the directional multiple testing procedure, which controls mdFDR. Guo et al. (2010) extend this procedure to the  comparison between many ordered treatments (or time points), where they concentrate mainly on the identification of the sign of the difference between mean values of the test statistic for consecutive treatments. According to their definition, the pattern for a given characteristic is wrongly identified if at least one of these signs is wrongly detected. Guo et al. (2010) prove that under this definition of the directional error their proposed procedure controlls mdFDR.

 While in principle the  multidimensional directional multiple testing problem can be addressed by the appropriate extension of the mixture model and application of  the Bayes or empirical Bayes approach, we are not aware of any specific Bayesian solution. Surely, this is a promising topic for a further research.

\section{Acknowledgement} We thank Professor Hira Koul for inviting us to write a paper for the Golden Jubilee volume of the Journal of the Indian Statistical Association. We also thank Professor Draper for sending us a copy of Draper and Krnjajic (2010). Warm thanks are due to two referees for penetrating and constructive criticism.

\end{document}